\documentclass[final,1p,times,authoryear]{elsarticle}
\usepackage{amsmath}
\usepackage{amsthm}
\usepackage{amssymb}
\usepackage{amsfonts}
\usepackage{xcolor} 
\usepackage[ruled,vlined]{algorithm2e} 

\SetKwInput{KwInput}{Input}
\SetKwInput{KwOutput}{Output}
\usepackage{listings} 
\newtheorem{theorem}{Theorem}

\newtheorem{definition}[theorem]{Definition}

\begin{document}

\begin{frontmatter}

\title{Unified computational approach to nilpotent algebra classification problems}

\author{Shirali Kadyrov}
\address{Suleyman Demirel University, Kaskelen, Kazakhstan}
\ead{shirali.kadyrov@sdu.edu.kz}

\author{Farukh Mashurov}
\address{Suleyman Demirel University, Kaskelen, Kazakhstan}
\address{Kazakh-British Technical University, Almaty, Kazakhstan}
\ead{farukh.mashurov@sdu.edu.kz}

\begin{abstract}
In this article, we provide an algorithm with Wolfram Mathematica code that gives a unified computational power in classification of finite dimensional nilpotent algebras using Skjelbred-Sund method. To illustrate the code, we obtain new finite dimensional Moufang algebras. 
\end{abstract}

\end{frontmatter}

\section{Introduction}

Let $\mathfrak A$ be a finite dimensional algebra over the filed of complex numbers $\mathbb C$ equipped with a bilinear product denoted by $x \cdot y$. Given a variety defined over certain polynomial identities, one of the classical research problems is to classify, up to isomorphism, all algebras within the variety with fixed dimension. Varieties of Associative, Alternative, Lie, Novikov, Jordan, Assosymmetric, Leibniz, Zinbiel, and Torkara algebras are some of the most studied varieties, see e.g.  \citep{ACK19,De18}, \citep{KKK19,GGKS19,KV19} and references therein. As described in the next section, classification involves several steps each requiring to symbolically solve systems of polynomial equations. Computational effort depends on many variables including the dimension of the variety, bilinear product, and number of identities and it may be required to solve system of more than $n^3$ polynomial equations with $n^2$ unknowns, where $n$ is the fixed dimension of the variety. Obviously, the computations could become cumbersome even when the dimension is as little as 2, see e.g. \citep{GR11}. As such, many authors rely on computer assisted classifications. 

Our goal in this article is to provide a unified algorithm with a code written in Wolfram Mathematica to ease the computation aspects of the classification problem for nilpotent algebras and illustrate with a new classification result. It is unified in the sense that a minimal effort is required to run the code, namely define the bilinear product and polynomial identities for the kind of algebraic variety being studied. Moreover, most of the computational steps are taken care of by the code as opposed to the previous works. Another significance of the current work is that it focuses on explaining the algorithms and the details of the full source code while the previous works were concentrated on algebraic classification with little attention to coding. The code is successfully tested on several previously obtained classification results. Moreover, it was used to classify assosymmetric algebras of dimension 4 and one generated assosymmetric algebras of dimension 5 and 6. We illustrate our code and obtain new 4-dimensional  Moufang  algebra. The central extensions of all two dimensional nilpotent algebras was recently obtained in \citep{CFK18}.

In the next section, we provide the background for steps needed to apply well-known Skjelbred-Sund classification method together with the algorithms that we follow in writing the code. In the follow up section \S~\ref{sec:ex} we provide new results to illustrate our unified symbolic computational approach. The original code by the authors is provided open access through \citep{Program}. Finally in \S~\ref{sec:conc} we conclude with possible future research directions.

\section{Skjelbred-Sund classification method}\label{sec:method}

Skjelbred-Sund classification method is one of the classical methods to classify finite dimensional nilpotent algebras which goes back to \cite{SS77}, where used central extensions of less than 6-dimensional Lie algebras to describe nilpotent  6-dimensional Lie algebras. More recently, \cite{HAC16} used Skjelbred-Sund in classification of all non-Lie central extensions of all 4-dimensional Malcev algebras and \cite{HA16} in classification of five-dimensional nilpotent Jordan algebras. For various results where the same method was used for algebraic classification of finite dimensional nilpotent algebras from different varieties we refer to \citep{HA18},  \citep{GKK19}, \citep{KGV19},  and references therein. 

For an excellent exposition of the method we refer to \citep{HA16} and \citep{HAC16}. Here, we simply review some important notions used in the classification.

We let $\mathfrak A$ denote an $n$ dimensional algebra  in certain variety $\mathfrak L$ defined by  set of polynomial identities and  $V$  be a vector space over $\mathbb{C}$. We define $Z^2_{\mathfrak L}(\mathfrak A,V)$ as the set of bilinear forms $\theta:\mathfrak A\ \times \mathfrak A  \longrightarrow {V}$  satisfying all identities in the variety $\mathfrak L.$ That is, $Z^2_{\mathfrak L}(\mathfrak A,V)$ is the set of closed bilinear forms, also known as \textit{cocycles}. For a linear map $f$ from $\mathfrak A$ to $V,$ we define a map $\delta f\colon \mathfrak A \times
\mathfrak A \longrightarrow V $ given by $\delta f(x,y)=f(xy)$. It is easy to see that $\delta f\in Z^2_{\mathfrak L}(\mathfrak A,V).$ Therefore, 
$${B^{2}_\mathfrak L}\left(
{\mathfrak A},{{ V}}\right) :=\left\{ \theta =\delta f\ | f\in {\rm Hom}\left( \mathfrak A,{{ V}}\right) \right\}$$
is a subspace of  $Z^2_{\mathfrak L}(\mathfrak A,V).$ The elements of ${B^{2}_\mathfrak L}\left(
{\mathfrak A},{{ V}}\right)$ are called \textit{coboundaries}. The second cohomolgy space is define as 
$${H^{2}_\mathfrak L}\left(\mathfrak A,V\right) := Z^2_{\mathfrak L}(\mathfrak A,V) \big/{B^{2}_\mathfrak L}\left(
{\mathfrak A},{{ V}}\right).$$

We let ${Aut}(\mathfrak A) $ denote the automorphism group of $\mathfrak A$ and take $\phi \in Aut(\mathfrak A).$ If we define an action of ${Aut}(\mathfrak A) $ on $Z^2_{\mathfrak L}(\mathfrak A,V)$ via $\phi\theta(x,y)=\theta(\phi(x),\phi(y))$ for every  $\theta \in Z^2_{\mathfrak L}(\mathfrak A,V),$ then necessarily $\phi\theta\in Z^2_{\mathfrak L}(\mathfrak A,V).$ That is, ${Aut}(\mathfrak A) $ acts on $Z^2_{\mathfrak L}(\mathfrak A,V).$ Furthermore, $\phi\theta\in B^{2}_\mathfrak L(
{\mathfrak A},{{ V}}) $ if and only if $\theta\in B^{2}_\mathfrak L(
{\mathfrak A},{{ V}}),$ so that ${Aut}(\mathfrak A) $ also acts on $H^{2}_\mathfrak L(
{\mathfrak A},{{ V}}).$

We let $V$ be a vector space of dimension $m$ and $\theta\in Z^2_{\mathfrak L}(\mathfrak A,V).$ With a multiplication ``$\cdot_{\mathfrak A_{\theta}}$'' on $\mathfrak A_{\theta}:=\mathfrak A\oplus V$ given by 
$$(x+v_1)\cdot_{\mathfrak A_{\theta}}(y+v_2)= x\cdot y + \theta(x,y) \text{ for all } x,y\in \mathfrak A \text{ and } v_1, v_2 \in V,$$
${\mathfrak A_{\theta}}$ becomes an algebra in the variety $\mathfrak L$, called \textit{m-dimensional central extension} of $\mathfrak A$ by vector space $V.$

Let $\{e_1,e_2,\dots,e_m\}$ be a basis for $V.$ Then, any $\theta\in Z^2_{\mathfrak L}(\mathfrak A,V)$ can be uniquely written as $\theta(x,y)=\sum^s_{i=1}\theta_i(x,y)e_i,$ for some $\theta_i\in Z^2_{\mathfrak L}(\mathfrak A,\mathbb{C}).$ The set ${Ann}(\theta):=\{ x\in {\mathfrak A}:\theta ( x, \mathfrak A )+ \theta (\mathfrak A ,x) =0\} $ is called
the {\it annihilator} of $\theta $. The {\it annihilator} of an  algebra ${\mathfrak A}$ is defined as
the ideal ${Ann}(  {\mathfrak A} ) =\{ x\in {\mathfrak A}:  x{\mathfrak A}+ {\mathfrak A}x =0\}$. 

We know that every finite-dimensional nilpotent algebra is a central extension of some nilpotent algebra of lower dimension,  see e.g. \citep{HAC16}. Therefore, to classify all nilpotent algebras of a fixed dimension in variety $\mathfrak L$, all we need is to classify cocycles of nilpotent algebras $\mathfrak A$ of lower dimension with condition $Ann(\mathfrak A)\cap Ann(\theta) = 0$ and central extensions that appear from them (see Lemma 5 in \cite{HAC16}).

Let $G_m(H^{2}_\mathfrak L(\mathfrak A,\mathbb C))$ be the set of all $m-$dimensional subspaces of $H^{2}_\mathfrak L(\mathfrak A,\mathbb C).$ We define $${T}_{m}({\mathfrak A}) =\{ {W}:=\langle [\theta _{1}] ,
[\theta _{2}] ,\dots,[ \theta _{m}]\rangle \in
{G}_{m}( { H^{2}_\mathfrak L}( {\mathfrak A},\mathbb C)) : \cap_{i=1}^{m}{\rm Ann}(\theta _{i})\cap{\rm Ann}({\mathfrak A}) = 0\}. $$

We define action of ${Aut}(\mathfrak A)$ on ${T}_{m}(\mathfrak A)$  via $\phi W=<[\phi\theta_1],[\phi\theta_2],\ldots,[\phi\theta_m]>$ for $\phi \in Aut(\mathfrak A)$ and $W \in {T}_{m}(\mathfrak A)$. The set $Orb(W)$ is the orbit of $W\in{T}_{m}(\mathfrak A)$ under automorphism group of $\mathfrak A.$

\begin{definition}
For $\mathfrak A=\mathfrak A_0 \oplus \mathbb C x $ where $x\in Ann(\mathfrak A),$ the subspace $\mathbb C x$ is called \emph{annihilator  component} of $\mathfrak A.$
\end{definition}

Assume that a vector space $V$ has dimension $m.$ For a given algebra $\mathfrak A$ in variety $\mathfrak L,$ we let $E(\mathfrak A,V)$ be the set of all  non-split central extensions of $\mathfrak A$ by $V,$ where a  non-split central extension is central extension of $\mathfrak A$ without annihilator components (see Definition 8 in \cite{KKP19}). 


\begin{theorem}[see in \cite{HAC16}, \cite{KKP19}]
There is a one-to-one correspondence between the set of $Aut(\mathfrak A)$-orbits on ${T}_{m}(\mathfrak A)$ and the set of isomorphism
classes of  $E(\mathfrak A,V).$
\end{theorem}

Finally, we have the following steps to construct from the algebra $\mathfrak A$ of dimension $n-m$ all non-split algebras in $\mathfrak L$ of dimension $n$ with $m$-dimensional annihilator. We need  to determine the following sets for given algebra $\mathfrak A:$
\begin{enumerate}
    \item Compute base for ${ Z^{2}_\mathfrak L}( \mathfrak A,\mathbb {C});$
     \item Compute base for ${ B^{2}_\mathfrak L}( \mathfrak A,\mathbb {C})$ and ${ H^{2}_\mathfrak L}( \mathfrak A,\mathbb {C}) ;$
  \item Compute ${Aut}({\mathfrak A});$
  \item Compute base for ${ Ann}(\mathfrak A)$ and ${ Ann}(\mathfrak A)\cap { Ann}(\theta);$
\item Compute ${\rm Aut}(\mathfrak A)$-orbits on ${T}_{m}({\mathfrak A}) ;$
\item  Construct the algebra in the variety $\mathfrak L$ associated with a
representative of each orbit.
\end{enumerate}

We now describe algorithms to handle steps from 1 to 4. The remaining two steps are work out by hand.

Let $\mathfrak A$ be an algebra with basis $\{e_i :i=1,2,\dots,n \}.$ 
We use the following notations: $\Delta_{i,j}$ is the bilinear form  $\Delta_{i,j}: \mathfrak A  \times \mathfrak A\longrightarrow \mathbb C$ such that \begin{equation}\label{kroneckerdelta}
    \Delta_{i,j}(e_l,e_k)=\delta_{il}\delta_{jk}.
\end{equation} The set $\{\Delta_{i,j}: 1\leq i,j\leq n\}$  is the basis of ${ Z^{2}_\mathfrak L}( \mathfrak A,\mathbb {C}).$ Every $\theta \in { Z^{2}_\mathfrak L}( \mathfrak A,\mathbb {C})$ can be uniquely written as $\theta=\sum^n_{1\leq i,j\leq n} \lambda_{i,j}\Delta_{i,j}$ where $\lambda_{i,j}\in \mathbb C.$

Now, we give algorithms to compute the above mentioned steps. The first algorithm shows how to compute ${Z^{2}_\mathfrak L}( \mathfrak A,\mathbb {C})$ given the dimension, the product rule, and the polynomial identities. It amounts to defining the symbolic equations and call the symbolic solver from the relevant programming language. 

\begin{algorithm}[H]\label{alg:Z2A}
\DontPrintSemicolon
  
  \KwInput{Dimension of your algebra: $n$, Bilinear product rule: $\text{pr}(\cdot,\cdot)$, Identities: $\{\text{Iden}_1,\dots,\text{Iden}_m\}$.}

    
    Define a symbolic basis for algebra $\mathfrak A$: $\{e_1,\dots,e_n\}$
    
    Define symbolic cocycles, a bilinear function: $\theta(\cdot,\cdot)=\sum_{i,j}^{n} \lambda_{i,j}\Delta_{i,j}(\cdot,\cdot)$ 
    Let $k$ be the number of variables used in the Identities.
    
    Define the system of symbolic $n^k m$ nonlinear equations:
    
     \For{$(i_1,\dots,i_k) \in \{1,\dots,n\}^k$}    
        { \For{j=1:m}
        {
        Define an $eq[i_1,\dots,i_k,j]$ that is obtained by applying $\theta$ to $\text{Iden}_j$ when $e[i_1], \dots, e[i_k]$ are substituted.
        }
        }
        
    Use built in `solve' function to solve the system $\{eq\}$.
    
 \KwOutput{Basis for ${ Z^{2}_\mathfrak L}( \mathfrak A,\mathbb {C})$: Z2A}
  
\caption{Algorithm to compute the basis for the ${ Z^{2}_\mathfrak L}( \mathfrak A,\mathbb {C})$}
\end{algorithm}

The next algorithm uses outcomes of Algorithm~\ref{alg:Z2A} together with the same inputs. In this case we aim to compute bases for ${ B^{2}_\mathfrak L}( \mathfrak A,\mathbb {C})$ and ${ H^{2}_\mathfrak L}( \mathfrak A,\mathbb {C})$. It does not require any tricks to obtain a basis for ${ B^{2}_\mathfrak L}( \mathfrak A,\mathbb {C})$ but simply write them down manually from the given polynomial identities. In terms of coding this means to ask the programming language to read coefficients of polynomial expressions. As for the second part, we recall that ${ B^{2}_\mathfrak L}( \mathfrak A,\mathbb {C}) \subset { Z^{2}_\mathfrak L}( \mathfrak A,\mathbb {C})$ and $ {H^{2}_\mathfrak L}\left(\mathfrak A,V\right) = Z^2_{\mathfrak L}(\mathfrak A,V) \big/{B^{2}_\mathfrak L}\left(
{\mathfrak A},{{ V}}\right)$. Thus, the problem of finding a basis for ${ H^{2}_\mathfrak L}( \mathfrak A,\mathbb {C})$ is equivalent to completing the basis of ${ Z^{2}_\mathfrak L}( \mathfrak A,\mathbb {C})$ given the basis of ${ B^{2}_\mathfrak L}( \mathfrak A,\mathbb {C})$.    

\begin{algorithm}[H]\label{alg:B2,H2}
\DontPrintSemicolon
  
  \KwInput{Dimension of your algebra: $n$, Bilinear product rule: $\text{pr}(\cdot,\cdot)$, Identities: $\{\text{Iden}_1,\dots,\text{Iden}_m\}$, Basis for $Z^2(\mathfrak A, \mathbb C): \{z_1,\dots,z_k\}$ }

From $\text{pr}(\cdot,\cdot)$ obtain a basis for $B^2(\mathfrak A):B2A=\{b_1,\dots,b_s\}$

Define empty set $H2A:=\{\};$

\For{i=1:k}{
  \If{$z_i \not\in \text{span}(B2A)$} {
   Add $z_i$ to $H2A$;
   
   Add $z_i$ to $B2A$;
   }
 }
\KwOutput{Basis for ${ B^{2}_\mathfrak L}( \mathfrak A,\mathbb {C})$ and ${ H^{2}_\mathfrak L}( \mathfrak A,\mathbb {C})$: B2A, H2A}
  
\caption{Algorithm to compute the bases for ${ B^{2}_\mathfrak L}( \mathfrak A,\mathbb {C})$ and ${ H^{2}_\mathfrak L}( \mathfrak A,\mathbb {C})$}
\end{algorithm}

 Computing the ${Aut}({\mathfrak A})$, Algorithm~\ref{alg:Aut}, is one of the main steps in the above described method and the one with large computational cost. We may represent an automorphism with an $n\times n$ invertible square matrix that respects the bilinear product rule. This requires to define symbolic matrix and define system of symbolic equations and finally call the solve function.

\begin{algorithm}[H]\label{alg:Aut}
\SetAlgoLined
\KwInput{Dimension of your algebra: $n$, Bilinear product rule: $\text{pr}(\cdot,\cdot)$}
 Define matrix GAut$_{n\times n}$;
 Define homomorphism function $F[\{x,y\}]=   pr[F[x],F[y]]-F[pr[x,y]]$\;
 Define mapping of basis by $F[e_i]=\sum^{n}_{j=1} \lambda_{i, j}e_j$\;
 \For{i,j; n}{Substitute basis to  $F[\{e_i,e_j\}]$ and define by Eq}
 Use "solve" to find 
  is set of solutions $\lambda_{i j}$ and define by Solution\;
 Obtained solutions substitute to matrix GAut, that is 
 
 \For{i=1; Length[Solutions]}{Mat[i]=GAut/.Solution[[i]]}
 
 Define set Automorphism=\{\};
 
 \If{Det[Mat[i]]!=0 }{Add to  Automorphism}
 
 \KwOutput{ Matrix forms of Automorphism group : Automorhism}
 
\caption{Algorithm to find the automorphism group}
\end{algorithm}

The next three algorithms are allocated for step 4 to compute annihilators. The Algorithm~\ref{alg:Action on Aut} computes the action of the automorphism group on $H^2_\mathfrak L(\mathfrak A,\mathbb C)$ and uses outcomes of Algorithm~\ref{alg:B2,H2} and Algorithm~\ref{alg:Aut}. Action of automorphism group defined by $\phi^T*M*\phi$ where $\phi\in Aut(\mathfrak A)$ and $M$ is matrix form of $H^2_\mathfrak L(\mathfrak A,\mathbb C).$

\begin{algorithm}[H]\label{alg:Action on Aut}
\SetAlgoLined
\KwInput{Automorphism group of algebra, basis of  $H^2_\mathfrak L(\mathfrak A,\mathbb C).$ }

 \For{i=1; Length[Automorphism]}{ActAut[i]=Transpose[Automorphism[[i]]].MatrixFormH2.Automorphism[[i]]}

 \KwOutput{ Action of automorphism group on $H^2_\mathfrak L(\mathfrak A,\mathbb C)$: ActAut[i]}
 
\caption{Algorithm of action of the automorphism group on $H^2_\mathfrak L(\mathfrak A,\mathbb C)$}
\end{algorithm}

Next algorithm uses Algorithm~\ref{alg:Action on Aut} to compute bases for ${ Ann}(\mathfrak A)$. Again one needs to define the system of polynomial equations and call the solver.

\begin{algorithm}[H]\label{alg:basis of Annihilator}
\SetAlgoLined
\KwInput{Multiplication of basis elements}

 Define linear combinations of basis elements of algebra by SumElemOFAnn;
 
 Define empty set by  ProdSumElemOFAnn;

 \For{i=1; dim}{
 
    pr[SumElemOFAnn, e[i]] == 0  add to ProdSumElemOFAnn;

  pr[e[i], SumElemOFAnn] == 0 add to ProdSumElemOFAnn;  }
 
 Solve ProdSumElemOFAnn after obtained solution put to SumElemOFAnn;

 \KwOutput{ Basis of Annihilator: SpanAnn.}
 
\caption{Finding basis of annihilator.}
\end{algorithm}

Finally, the last algorithm below uses outcome of Algorithm \ref{alg:basis of Annihilator} and gives conditions of ${ Ann}(\mathfrak A)\cap { Ann}(\theta)=0:$

\begin{algorithm}[H]\label{alg:intersection of Ann and theta}
\SetAlgoLined
\KwInput{$\mathfrak A$ is $n$ dimensional algebra,  basis  of ${ Ann}(\mathfrak A)$ and basis of $H^2_\mathfrak L(\mathfrak A,\mathbb C)$ }
Define the following :

Linear combinations of basis elements of ${ Ann}(\mathfrak A)$ by SpanAnn1.

Define $\Delta_{i,j}$ as (\ref{kroneckerdelta});

 OpenBraket[x\_, y\_] := x[y];
 
 ConditionOfAnnAndTheta = \{\};
 
 \For{j;n+1}{\For{i;Length[H2]}{ Sum[$\alpha_i$OpenBraket[H2[[i]]],\{SpanAnn1[[1]], e[j]\}]==0 add to ConditionOfAnnAndTheta;
 
 Sum[$\alpha_i$OpenBraket[H2[[i]]],\{e[j],SpanAnn1[[1]]\}]==0 add to ConditionOfAnnAndTheta;}}
 
\KwOutput{Condition of ${Ann}(\mathfrak A) \cap{ Ann}(\theta)=0:$ ConditionOfAnnAndTheta .}
 
\caption{Intersection condition of ${ Ann}(\mathfrak A)$ and ${ Ann}(\theta):$ }
\end{algorithm}

\section{Application: classification of four dimensional nilpotent Moufang algebras.} \label{sec:ex}

To illustrate the code we now obtain new 4-dimensional algebra for the variety of Moufang algerbas.  For the detailed previous study on Moufang  algebras see in \cite{L93}, \cite{ShP04} and here we adopt the definition given by \cite{Loday}.  

An algebra $\mathfrak M$ is called Moufang algebra if it satisfies the following polynomial identities:
\begin{equation}\label{Moufangident}
\begin{array}{c}
(x,y,z)=-(z,y,x),\\
((xy)z)t + ((zy)x)t = x(y(zt)) + z(y(xt)),\\
t(x(yz) + z(yx)) = ((tx)y)z + ((tz)y)x,\\
(xy)(tz) + (zy)(tx) = (x(yt))z + (z(yt))x,
\end{array}\end{equation}
where $(x,y,z)=(xy)z-x(yz).$

Firstly, we need to give for the code some information about Moufang algebra in the following form: 
     
     We first input the number of identities denoted by $n.$ In our case we set 
    \begin{lstlisting}
    n = 4;
    \end{lstlisting}
    
    The multiplication of basis elements defined by $\text{pr}[*,*]$. The identities defined such functions $\text{Ident}[i][\{x\_,y\_,z\_,...\}]$, where $i \in \{1,..,n\}$ and first product should be written with capital ``Pr''. Given identities should be in homogeneous multilinear form. Also, we denote the length of monomials in identities by LengthOfMonomial[i] where $i \in \{1,..,n\}.$ Since, we consider Moufang algebra we input the identities (\ref{Moufangident}) in the following form
      
\begin{lstlisting}
Ident[1][{x_,y_,z_}]:= 
 Pr[pr[x,y],z]-Pr[x,pr[y,z]]+Pr[pr[z,y],x]-Pr[z,pr[y,x]]
LengthOfMonomial[1] = 3;

Ident[2][{x_,y_,z_,t_}]:= 
    Pr[pr[pr[x,y],z],t]-Pr[x,pr[y,pr[z,t]]]+
        Pr[pr[pr[z,y],x],t]-Pr[z,pr[y,pr[x,t]]]
LengthOfMonomial[2] = 4;

Ident[3][{x_,y_,z_,t_}]:= 
    Pr[pr[pr[t,x],y],z]-Pr[t,pr[x,pr[y,z]]]+ 
        Pr[pr[pr[t,z],y],x]-Pr[t,pr[z,pr[y,x]]]
LengthOfMonomial[3] = 4;

Ident[4][{x_,y_,z_,t_}]:= 
    Pr[pr[x,y],pr[t,z]]+Pr[pr[z,y],pr[t,x]]-
        Pr[pr[x,pr[y,t]],z]-Pr[pr[z,pr[y,t]],x]
LengthOfMonomial[4] = 4;
 
        \end{lstlisting}

\subsection{Choosing three dimensional Moufang algebras from the literature}
As mentioned before, various three dimensional nilpotent algebra classifications were considered in \citep{CFK18} including Moufang algebras. For our purposes, we need to extract the three dimensional algebras from this article.

To this end, we need to input dimension of algebra, in the code it is given as $dim.$ So,
        \begin{lstlisting}
        dim =3;
        \end{lstlisting}
 Next step, we input multiplication table of algebra in the form 
   $\text{pr}[e[i],e[j]]=e[k]$ where $i,j\in \{1,..,n\}.$ If for some $i,j$ the product ${\rm{pr}}[e[i],e[j]]=0$ then drop it. For instance,
   
   \begin{lstlisting}
        pr[e[1],e[1]]:=e[2];
        pr[e[1],e[2]]:=e[3];
        pr[e[2],e[1]]:=e[3];
\end{lstlisting} 

Now the code has all the information that it needs. Next we select all three dimensional nilpotent  Moufang algebras  from \citep{CFK18} by using the following code: 

\begin{lstlisting}
ChangeTopr[x_]:= x/.{Pr->pr}
Unident={};
Do[Unident= 
  Union[Unident, 
   Expand[ChangeTopr[
     Map[Ident[i],Tuples[basis,LengthOfMonomial[i]]]]]],{i,1,n}]
If[Unident == {0}, 
 Print[Grid[{Text@
      Style[#,"TableHeader"]&/@{"The identities hold true."}}, 
   Frame -> All]], 
 Print[Grid[{Text@
      Style[#,"TableHeader"]&/@{"The identities do not hold."}},
   Frame -> All]]]
\end{lstlisting}
 As result we have the following 3-dimensional nilpotent Moufang algebras: 

$$\begin{array}{ll llllllllllll}
{\mathfrak M}^3_{01} &:& e_1 e_1 = e_2,  & e_1 e_2=e_3,  & e_2 e_1 = e_3;\\
{\mathfrak M}^{3}_{02} &:& e_1 e_1 = e_2;\\
{\mathfrak M}^{3}_{03} &:& e_1 e_2 = e_3,  & e_2 e_1=e_3;  \\ 
{\mathfrak M}^{3}_{04} &:& e_1 e_2=e_3, & e_2 e_1=-e_3;   \\
{\mathfrak M}^{3}_{05}(\lambda) &:& e_1 e_1=\lambda e_3, & e_2 e_1=e_3 & e_2e_2=e_3;  \\

\end{array} $$

\subsection{Second cohomology of three dimensional nilpotent Moufang algebras.}

Here, we apply our code constructed by Algorithm \ref{alg:Z2A} and Algorithm \ref{alg:B2,H2} to obtain ${ Z^{2}_\mathfrak L}( \mathfrak A,\mathbb {C})$, ${ B^{2}_\mathfrak L}( \mathfrak A,\mathbb {C})$, and ${ H^{2}_\mathfrak L}( \mathfrak A,\mathbb {C}).$

{\tiny
$$
\begin{array}{|l|l|l|l|}
\hline
\mathfrak M  & {Z^{2}}\left( {\mathfrak M}, \mathbb C\right)  & { B}^2({\mathfrak M}, \mathbb C) & { H}^2({\mathfrak M}, \mathbb C) \\
\hline

{\mathfrak M}^{3}_{01} &  
\Big\langle\begin{array}{l}\Delta _{1,1},\Delta _{1,2}+\Delta _{2,1},\Delta _{1,3}+\Delta _{2,2}+\Delta _{3,1} \end{array}\Big\rangle& 
\Big\langle \begin{array}{l}\Delta _{1,1},\Delta _{1,2}+\Delta _{2,1} \end{array}\Big\rangle&
\Big\langle \begin{array}{l}[\Delta _{1,3}]+[\Delta _{2,2}]+[\Delta _{3,1}] \end{array}\Big\rangle\\
\hline
{\mathfrak M}^{3}_{02}  & \Big\langle\begin{array}{l}\Delta _{1,1},\Delta _{1,3},\Delta _{1,2}+\Delta _{2,1},\Delta _{3,1},\Delta_{2,3}+\Delta _{3,2},\Delta _{3,3}\end{array} \Big\rangle & \Big\langle \begin{array}{c}
     \Delta _{11}
\end{array} \Big\rangle& \Big\langle\begin{array}{c}
 [\Delta _{1,3}],[\Delta _{1,2}]+[\Delta _{2,1}],[\Delta _{3,1}],[\Delta _{2,3}]+[\Delta _{3,2}],[\Delta _{3,3}] 
\end{array}\Big\rangle\\
\hline
{\mathfrak M}^{3}_{03} &  
\Big\langle\begin{array}{l} \Delta _{1,1},\Delta _{1,2},\Delta _{2,1},\Delta _{2,2},\Delta _{1,3}+\Delta _{3,1},\Delta _{2,3}+\Delta _{3,2}\end{array}\Big\rangle& \Big\langle\begin{array}{c}
    
\Delta _{1,2}+\Delta _{2,1}\end{array} \Big\rangle& \Big\langle\begin{array}{c}
    [\Delta _{1,1}],[\Delta _{2,1}],[\Delta _{2,2}],[\Delta _{1,3}]+[\Delta _{3,1}],[\Delta _{2,3}]+[\Delta _{3,2}]\end{array}\Big\rangle\\
\hline

{\mathfrak M}^{3}_{04}&\Big\langle\begin{array}{l} \Delta _{1,1},\Delta _{1,2},\Delta _{2,1},\Delta _{2,2},\Delta _{1,3}-\Delta _{3,1},\Delta _{2,3}-\Delta _{3,2} \end{array}\Big\rangle & \Big\langle\begin{array}{l}\Delta _{1,2}-\Delta _{2,1} \end{array}\Big\rangle& \Big\langle\begin{array}{l}[\Delta _{1,1}],[\Delta _{1,2}],[\Delta _{2,2}],[\Delta _{1,3}]-[\Delta _{3,1}],[\Delta _{2,3}]-[\Delta _{3,2}]\end{array}\Big\rangle\\
\hline

{\mathfrak M}^{3}_{05}(\alpha) & \Big\langle\begin{array}{l} \Delta _{1,1},\Delta _{1,2},\Delta _{2,1},\Delta _{2,2} \end{array}\Big\rangle & \Big\langle\begin{array}{l}\lambda\Delta _{1,1}+\Delta _{2,1}+\Delta _{2,2} \end{array}\Big\rangle& \Big\langle\begin{array}{l}[\Delta _{1,1}],[\Delta _{1,2}],[\Delta _{2,2}]\end{array}\Big\rangle\\
\hline
\end{array}$$}

We see from the above table that the second cohomology spaces and automorphism groups of   ${\mathfrak M}^{3}_{01},{\mathfrak M}^{3}_{02},{\mathfrak M}^{3}_{03},{\mathfrak M}^{3}_{04},{\mathfrak M}^{3}_{05}(\lambda)$ algebras coincide with  $\mathfrak J ^3_{01},\mathfrak J ^{3*}_{01},\mathfrak J ^{3*}_{02},\mathfrak J ^{3*}_{03},\mathfrak J ^{3*}_{04}$ algebras in \citep{{JKKh19}} respectively. As such, it is clear that these algebras have the same central extensions. In other words, our calculations are expected to produce the same four dimensional algebras. Indeed, that is what we will get next which in particular gives us an opportunity to verify our code.

\subsection{Computation of central extensions of ${\mathfrak M}^{3}_{01}$. }

The next part is computing ${\rm Aut} ({\mathfrak M}^{3}_{01}).$ This part constructed by Algorithm \ref{alg:Aut} and we get the automorphism group ${\rm Aut} ({\mathfrak M}^{3}_{01})$ denoted by the set ``Automorphism'' in the code. In our partucular case, we get the following output:
$$\left(
\begin{array}{ccc}
 \lambda_{1,1} & 0 & 0 \\
 \lambda_{2,1} & \lambda_{1,1}^2 & 0 \\
 \lambda_{3,1} & 2 \lambda_{1,1} \lambda_{2,1} & \lambda_{1,1}^3 \\
\end{array}
\right),$$
By Algorithm \ref{alg:Action on Aut} the action of $\phi\in{\rm Aut} ({\mathfrak M}^{3}_{01})$ on the subspace 
$[\theta]= \alpha_1 ([\Delta _{1,3}]+[\Delta _{2,2}]+[\Delta _{3,1}]) $ gives the following table:
{\tiny
$$
\begin{array}{|l|l|l|l|l|}
\hline
{\mathfrak M} & {\rm Aut} ({\mathfrak M}) &\text{ Matrix form of } { H}^2({\mathfrak M}, \mathbb C)  & \text{Action of automorphism group on  } { H}^2({\mathfrak M}, \mathbb C) & \text{Where $\alpha^*_i $'s are equal:}\\
\hline
{\mathfrak M}^{3}_{01} &
\left(
\begin{array}{ccc}
 \lambda_{1,1} & 0 & 0 \\
 \lambda_{2,1} & \lambda_{1,1}^2 & 0 \\
 \lambda_{3,1} & 2 \lambda_{1,1} \lambda_{2,1} & \lambda_{1,1}^3 \\
\end{array}
\right) &  

\left(
\begin{array}{ccc}
 0 & 0 & \alpha_1 \\
 0 & \alpha_1  & 0 \\
 \alpha_1 & 0 & 0 \\
\end{array}
\right) &

 \left(
\begin{array}{ccc}
 \beta_1 & \beta_2 & \alpha ^*_1 \\
 \beta_2 & \alpha ^*_1 & 0 \\
 \alpha ^*_1 & 0 & 0 \\
\end{array}
\right) & \alpha ^*_1=\alpha_1 \lambda_{1,1}^4 \\\hline

\hline
\end{array}$$
}

Also, one of the steps of the construction a non-split algebra is computing basis of annihilator and conditions required for $Ann([\theta])\cap Ann({\mathfrak M}^{3}_{01})=0.$ By Algorithm~\ref{alg:basis of Annihilator} 
our code gives output SpanAnn which is the set of basis elements for $Ann({\mathfrak M}^{3}_{01})$. And by Algorithm \ref{alg:intersection of Ann and theta} our code computes the conditions for $Ann([\theta])\cap Ann({\mathfrak M}^{3}_{01})=0.$
That is, any element in $u\in Ann({\mathfrak M}^{3}_{01})$ can be expressed by $\lambda_3e_3$ where $\lambda_3\in \mathbb {C}$ and it is defined in the code by SpanAnn1. Also, every $[\theta]\in {H}^2({\mathfrak M} \mathbb C)$ is linear combinations of $\alpha_1 ([\Delta _{1,3}]+[\Delta _{2,2}]+[\Delta _{3,1}]).$ Note that $\Delta _{i,j}$ is function defined as (\ref{kroneckerdelta}). We just check if $u=\lambda_3e_3\in $ SpanAnn1, when the following conditions $\theta(e_j,u)=0$ and $\theta( u,e_j)=0$ hold, where $j\in \{1,\ldots n\}.$ So, we get the following output:

$$
\begin{array}{|l|l|}
\hline
Ann({\mathfrak M}^{3}_{01})  &\text{Condition of } Ann([\theta])\cap Ann({\mathfrak M}^{3}_{01})\\
\hline
\{e_3 \} & \alpha _1 \lambda _3=0\\
\hline
\end{array}$$

All of the above calculations are performed within seconds by using our code. At this stage, we need to obtain different orbits by hand. 
Since, we are interested in $Ann([\theta])\cap Ann({\mathfrak M}^{3}_{01})=0$  we see that from above table that this holds whenever  $\alpha_1\neq0.$ If we take $\lambda_{1,1}=\frac{1}{\sqrt[4]{\alpha_1}}$ we get $\langle[\Delta _{1,3}]+[\Delta _{2,2}]+[\Delta _{3,1}]\rangle$, which finally leads to new 4-dimensional nilpotent Moufang algebra from ${\mathfrak M}^{3}_{01}:$
$$\begin{array}{|l|llllll|}\hline
   {\mathfrak M}^{4}_{01}  &  e_1e_1=e_2 & e_1e_2=e_3 & e_2e_1=e_3 & e_1e_3=e_4 & e_2e_2=e_4 & e_3e_1=e_4\\ \hline
\end{array}$$


\section{Conclusion} \label{sec:conc}

In this article, we explained a unified approach to produce new nilpotent algebras which is one of the active research area in algebra. In part, we use the function ``solve'', the built in symbolic solver of Wolfram Mathematica, in handling the system of polynomial equations. This is the main function that takes most of the compilation time. The codes written by the authors in other software including Matlab and Python gave worse results in terms of the running time and in some cases failing to provide any solutions. It is no doubt that solving system of symbolic nonlinear equations is not an easy task even for such advanced programming language. However, one may dig into the ``solve'' function and may improve to provide faster results for the kind of problems considered here.

\bibliographystyle{elsarticle-harv}

\end{document}